\documentclass[a4paper,11pt]{article}
\usepackage[latin1]{inputenc}
\usepackage{graphicx}
\usepackage{amsmath}
\usepackage{amsfonts}
\usepackage{here}
\newtheorem{thm}{Theorem}

\newtheorem{cor}[thm]{Corollary}

\newcommand{\white}{\qquad{\framebox{\rule{0pt}{4pt}}}}
\newcommand{\owari}{\hfill\white}

\def\va{\textbf{\textit{a}}}
\def\vb{\textbf{\textit{b}}}
\def\vvarepsilon{\boldsymbol{\varepsilon}}

\begin{document}


\begin{center}
{\Large A Correspondence between Chord Diagrams and Families of $0$-$1$ Young Diagrams}
\end{center}
\begin{center}
{\large Tomoki Nakamigawa\footnote{This work was supported by JSPS KAKENHI Grant Number 19K03607.
The author states that there is no conflict of interest.
All data generated or analyzed during this study are included in this published article.
\\Email: {\texttt{nakami@info.shonan-it.ac.jp}}}
}
\end{center}
\begin{center}
{
Department of Information Science \\%
Shonan Institute of Technology \\%
1-1-25 Tsujido-Nishikaigan, Fujisawa, \\%
Kanagawa 251-8511, Japan \\%
e-mail: {nakami@info.shonan-it.ac.jp}
}
\end{center}
\begin{abstract}
A chord diagram is a set of chords in which no pair of chords has a common endvertex.
For a chord diagram $E$ having a crossing $S = \{ ac, bd \}$, by the chord expansion of $E$ with respect to $S$, we have two chord diagrams $E_1 = (E\setminus S) \cup \{ ab, cd \}$ and $E_2 = (E\setminus S) \cup \{ da, bc \}$.
Starting from a chord diagram $E$, by iterating expansions, we have a binary tree $T$ such that $E$ is a root of $T$ and a multiset of nonintersecting chord diagrams appear in the set of leaves of $T$.
The number of leaves, which is not depending on the choice of expansions, is called the chord expansion number of $E$.
A $0$-$1$ Young diagram is a Young diagram having a value of $0$ or $1$ for all boxes.
This paper shows that the chord expansion number of some type counts the number of $0$-$1$ Young diagrams under some conditions.
In particular,  it is shown that the chord expansion number of an $n$-crossing, which corresponds to the Euler number, equals the number of $0$-$1$ Young diagrams of shape $(n,n-1,\ldots,1)$ such that each column has at most one $1$ and each row has an even number of $1$'s.
\end{abstract}
keywords: chord diagram, chord expansion, Young diagram, Euler number, alternating permutation, increasing tree
\medskip

\section{Introduction}
A {\it chord} is a line segment whose endvertices are lying on a circle.
Two chords are called {\it independent} if they have no common endvertex. 
A {\it chord diagram} is a set of mutually independent chords.
A {\it crossing} is a pair of chords intersecting each other within the circle.
A chord is called {\it isolated} if it intersects no other chord.
A chord diagram is called {\it nonintersecting} if it has no crossing.
Let $E$ be a chord diagram having a crossing $S=\{ ac, bd \}$.
The {\it chord expansion} of $E$ with respect to $S$ consists in replacing $E$ with two chord diagrams $E_1 = (E \setminus S) \cup \{ ab, cd \}$ and  $E_2 = (E \setminus S) \cup \{ da, bc \}$.
Note that the number of crossings of $E_i$ is less than that of $E$ for $i=1,2$.
Hence, by iterating chord expansions, we have eventually a multiset of nonintersecting chord diagrams.
It is known that the multiset is uniquely determined not depending on the choice of expansions
 \cite{N2016}. 
We call the multiset the {\it chord expansion set} of $E$ and denote it by ${\cal NCD}(E)$.

A chord diagram is a combinatorial tool to study the representation theory of symmetric groups \cite{HJO2022, KR1984, Rhoades2019, RT2019}.
The {\it Specht module} $S^{(n,n)}$, one of the representations for the symmetric group ${\cal S}_{2n}$, has two important bases, the {\it polytabloid(Specht) basis}, and the {\it web basis}.
For a chord diagram $E$ having $n$ chords, let $[2n]=\{ 1,2,\ldots, 2n \}$ be the vertex set of $E$, where $1,2,\ldots,2n$ is arranged in anticlockwise order along the circumference.
Then $E$ is considered as a (perfect) {\it matching} on $[2n]$.
A matching is called {\it nonnesting} if there is no pair of chords $ij$ and $k\ell$ such that $k < i < j < \ell$.
The set of nonnesting matchings ${\rm NN}_{2n}$ is considered as an index set for the polytabloid basis $\{ v_E : E \in {\rm NN}_{2n} \}$.
On the other hand, the set of nonintersecting matchings ${\rm NC}_{2n}$ is an index set for the web basis $\{ w_F : F \in {\rm NC}_{2n} \}$.
Rhoades \cite{Rhoades2019} showed that the entry $c_{EF}$ of the transition matrix from the polytabloid basis $\{ v_E \}$ to the web basis $\{ w_F \}$ equals the number of occurrences of $F$ in ${\cal NCD}(E)$.
Hwang, Jang, and Oh \cite{HJO2022} introduced a new class of permutation, a web permutation, and showed that a set of web permutations on $[n]$ corresponds to ${\cal NCD}(E)$ for $E \in {\rm NN}_{2n}$.

\begin{figure}[H]
\begin{center}
\includegraphics[scale=0.3]{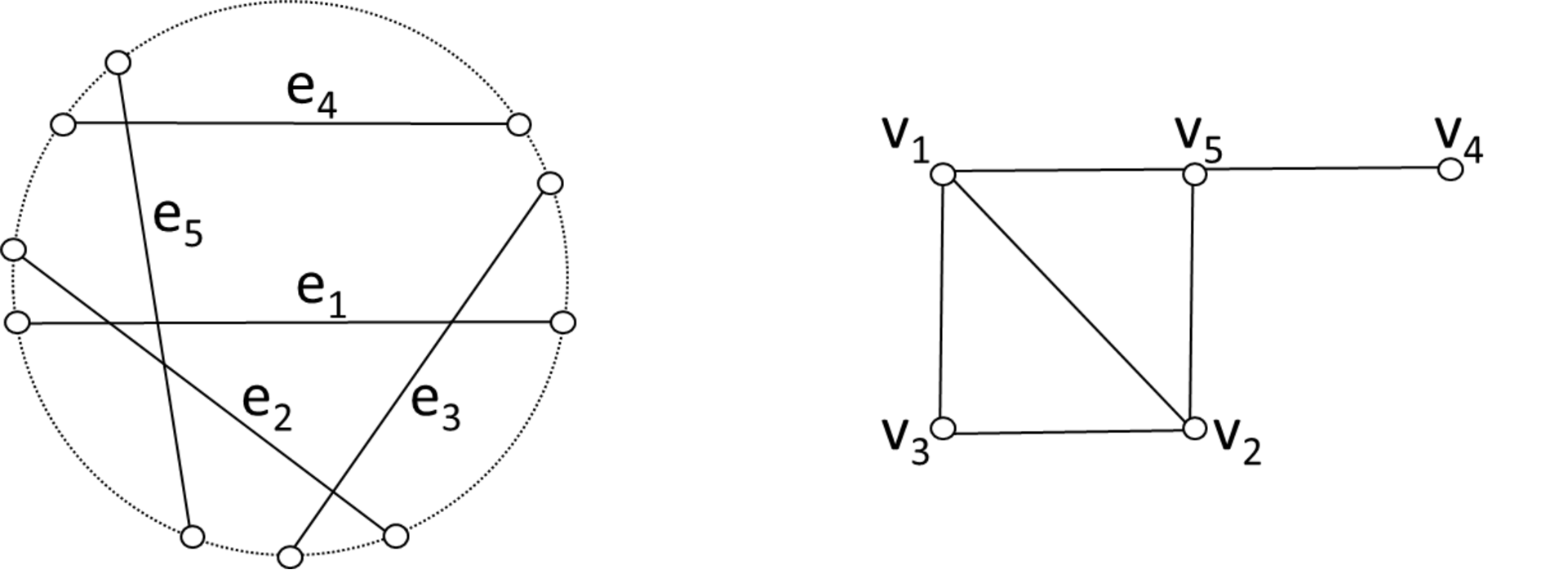}\\
\caption{A chord diagram $E$(left) and its associated intersection graph ${\rm Int}(E)$(right), where $v_i$ coressponds to $e_i$ for $1 \le i \le 5$.}\label{fig2311_intgraph}
\end{center}
\end{figure}

The cardinality of ${\cal NCD}(E)$ as a multiset is called the {\it chord expansion number} of $E$ and is denoted by $ex(E)$.
For a chord diagram $E$, let ${\rm Int}(E)$ be the {\it intersection graph} of $E$, namely the vertex set of ${\rm Int}(E)$ is $E$ and an edge of ${\rm Int}(E)$ corresponds to a crossing of $E$(see Fig. \ref{fig2311_intgraph}).
Chord diagrams $E_1$ and $E_2$ are called {\it isomorphic}, denoted by $E_1 \cong E_2$, if ${\rm Int}(E_1)$ and ${\rm Int}(E_2)$ are isomorphic as graphs.

The following results are known about the chord expansion number:
\begin{itemize}
\item For a chord diagram $E$, $ex(E)$ is equal to $T ({\rm Int}(E); 2, -1)$, where $T(G; x, y)$ is the {\it Tutte polynomial} of a graph $G$ \cite{NS2018}. 
\item For an $n$-{\it crossing} $C_n$, a set of mutually intersecting $n$ chords, ${\rm Int}(C_n)$ is a complete graph $K_n$ with $n$ vertices. 
The {\it Euler number} $Eul_n$, where $(Eul_n)_{n\ge 0} = (1,1,1,2,5,16,61,272,\ldots )$, counts the number of {\it down-up permutations} of $[n]$.
For $n \ge 1$, we have $ex(C_n) = Eul_{n+1}$ \cite{M2008, N2016}.
\item For two nonnegative integers $k$ and $n$ with $k \le n$, let $E_{n,k}$ be a chord diagram with $n+1$ chords, in which there is an $n$-crossing $C_n$ with an extra chord $e$ such that $e$ crosses exactly $k$ chords of $C_n$. 
For $1 \le k \le n$, the {\it Entringer number} $Ent_{n,k}$ is defined as the number of down-up permutations of $[n+1]$ with the first term $k+1$.
For $0 \le k \le n$, we have $ex(E_{n,k}) = Ent_{n+2,k+1}$ \cite{N2016}.
\item For an $(m,n)$-{\it crossing} $C_{m, n}$, a set of $m$ parallel chords $F_1$ and a set of $n$ parallel chords $F_2$ and $e_1$ and $e_2$ crosses for all $e_1 \in F_1$ and $e_2 \in F_2$, ${\rm Int}(C_{m, n})$ is a complete bipartite graph $K_{m, n}$.
Let $a_{m, n} = T(K_{m, n}; 2,-1)$.
Then its generating function is\\ $\sum_{m \ge 0} \sum_{n \ge 0} a_{m, n} \frac{x^m y^n}{m! n!} = 1/(\cosh x \cosh y - \sinh x \sinh y)$ \cite{GMMN2013}.
\end{itemize}
For more information on $T(G; 2,-1)$, see \cite{GMMN2013}.  
For other results on the chord expansion, see \cite{N2016, NS2018, N2020}.

In this paper, we focus on the enumerative aspect of chord expansion, and we study a relation between the chord expansion number and $0$-$1$ Young diagrams.
For an integer partition $\va = (a_1, a_2, \ldots, a_m)$ with $a_1 \ge a_2 \ge \cdots \ge a_m$, {\it Young diagram} of shape $\va$ is a set of boxes arranged in a matrix manner, in which the $i$th row has $a_i$ boxes for $1 \le i \le m$. 
A $0$-$1$ {\it Young diagram} is a Young diagram in which each box has a value of $0$ or $1$.


\section{Main Results}
\begin{figure}
\begin{center}
\includegraphics[scale=0.3]{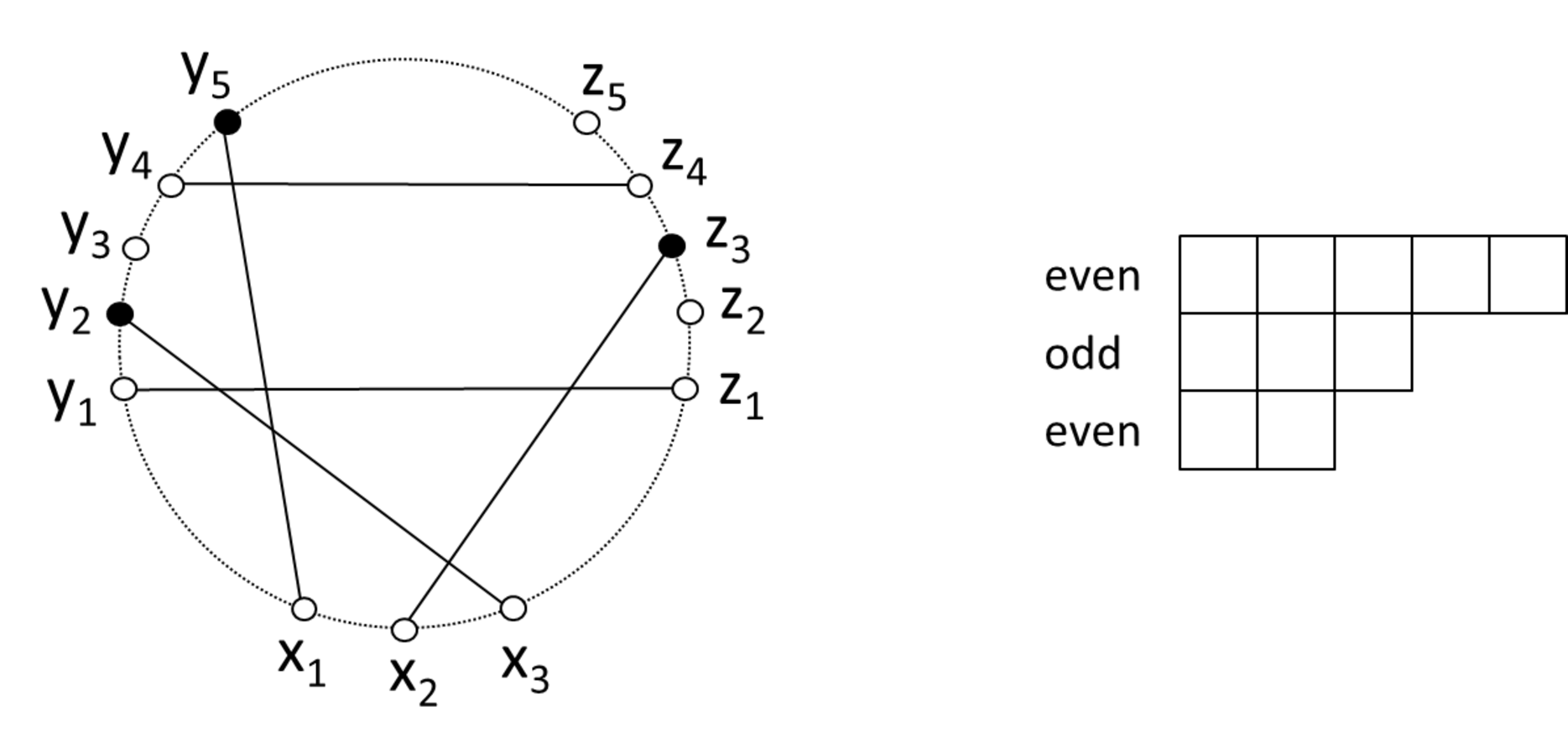}\\
\caption{A chord diagram $E((5,3,2), (0,1,0))$(left) and a Young diagram of shape $(5,3,2)$(right).}\label{fig2301_E532010}
\end{center}
\end{figure}

For a pair of positive integers $m$ and $n$ with $m \le n$, let $\va = (a_1, a_2, \ldots, a_m)$ be an integer partition, where $n = a_1 > a_2 > \cdots > a_m \ge 1$, and let $\vvarepsilon = (\varepsilon_1, \varepsilon_2,  \ldots, \varepsilon_m)$ be a $0$-$1$ sequence of length $m$.
Young diagram of shape $\va$ has $m$ rows and $n$ columns. 
Let us define a chord diagram $E(\va, \vvarepsilon)$ corresponding to the above-mentioned Young diagram as follows.
Endvertices of chords of $E(\va, \vvarepsilon)$ are contained in $X \cup Y \cup Z$, where $X = \{ x_1, x_2, \ldots, x_m \}$, $Y=\{ y_1, y_2, \ldots, y_n \}$, $Z=\{ z_1, z_2, \ldots, z_n \}$, and
$x_1, x_2, \ldots, x_m, $
$z_1, z_2, \ldots, z_n, $
$y_n, y_{n-1}, \ldots, y_1$ are located on the circumference in anticlockwise order.
$E(\va, \vvarepsilon)$ is defined as 
$\{ x_i y_j : 1 \le i \le m, a_i = j, \varepsilon_i = 0 \} \cup \{ x_i z_j : 1 \le i \le m, a_i = j, \varepsilon_i = 1 \} \cup \{ y_j z_j : j \in \{ 1,2, \ldots, n \} \setminus \{ a_1, \ldots, a_m \} \} $.
(See Fig. \ref{fig2301_E532010}.) 
Note that in the definition of $E(\va, \vvarepsilon)$, a set of $m$ isolated vertices always exists.
We eliminate these isolated vertices and consider $E(\va, \vvarepsilon)$ as a chord diagram having $n$ chords.

\begin{figure}[H]
\begin{center}
\includegraphics[scale=0.3]{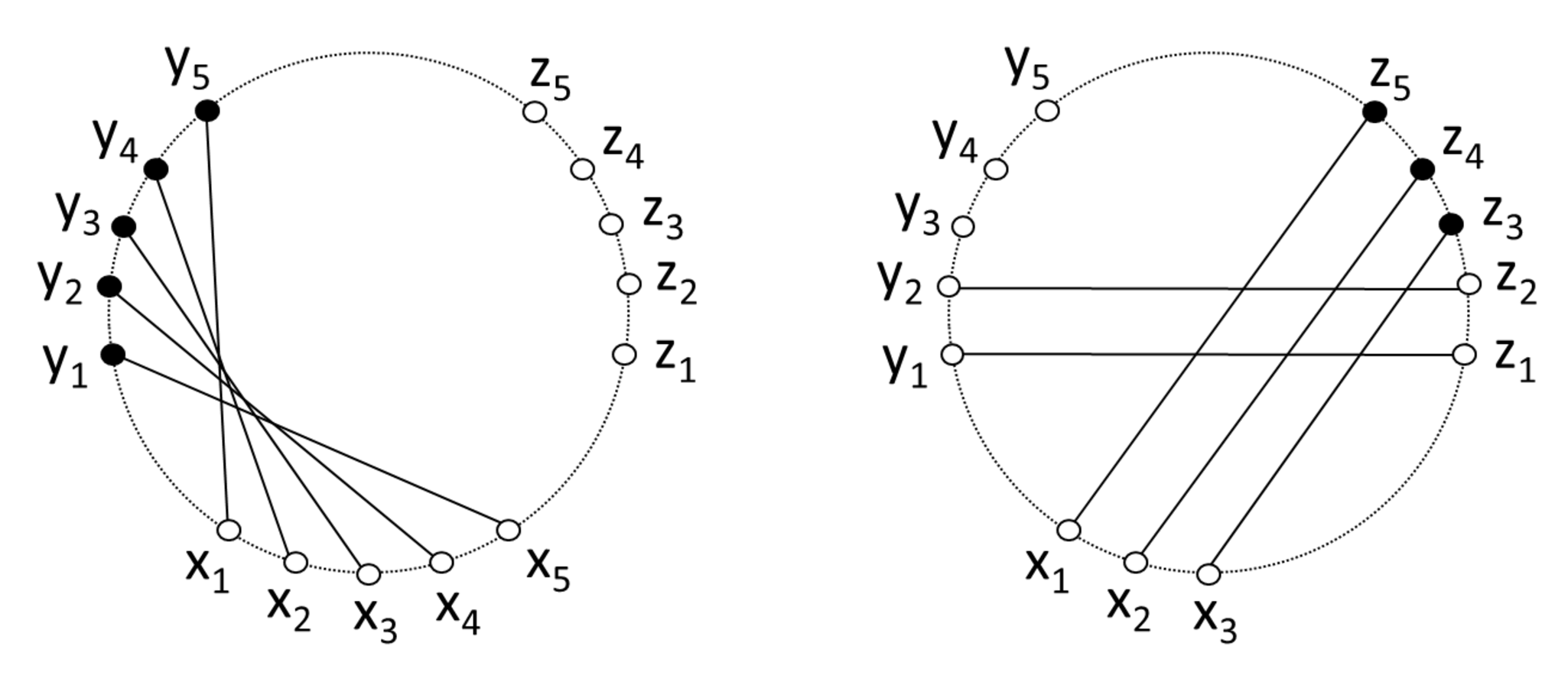}\\
\caption{$E((5,4,3,2,1), (0^5))$ (left) and $E((5,4,3), (1^3))$ (right).}\label{fig2311_n_cross_mn_cross}
\end{center}
\end{figure}

Let ${\cal E}_0$ be the family of chord diagrams represented by $E(\va, \vvarepsilon)$ for some pair of $\va$ and $\vvarepsilon$.
Note that ${\cal E}_0$ does not contain all chord diagrams.
For example, if $E \in {\cal E}_0$ contains two $2$-crossings $S$ and $T$ with $S \cap T = \emptyset$, then there exists a $2$-crossing $\{ e, f \}$ with $e \in S$ and $f \in T$. 
Nevertheless, ${\cal E}_0$ contains some important families of chord diagrams.
For a $0$-$1$ sequence, if a symbol $\varepsilon$ repeats $t$ times, we simply write $(\varepsilon^t)$.
For $\va = (n, n-1, \ldots, 1)$ and $\vvarepsilon = (0^n)$, we have $E(\va,\vvarepsilon)$ is an $n$-crossing $C_n$ (see Fig. \ref{fig2311_n_cross_mn_cross}(left)).
For $\va = (n+m, n+m-1, \ldots, n+1)$ and  $\vvarepsilon = (1^m)$, we have $E(\va,\vvarepsilon)$ is an $(m,n)$-crossing $C_{m,n}$ (see Fig. \ref{fig2311_n_cross_mn_cross}(right)).

\begin{figure}[H]
\begin{center}
\includegraphics[scale=0.3]{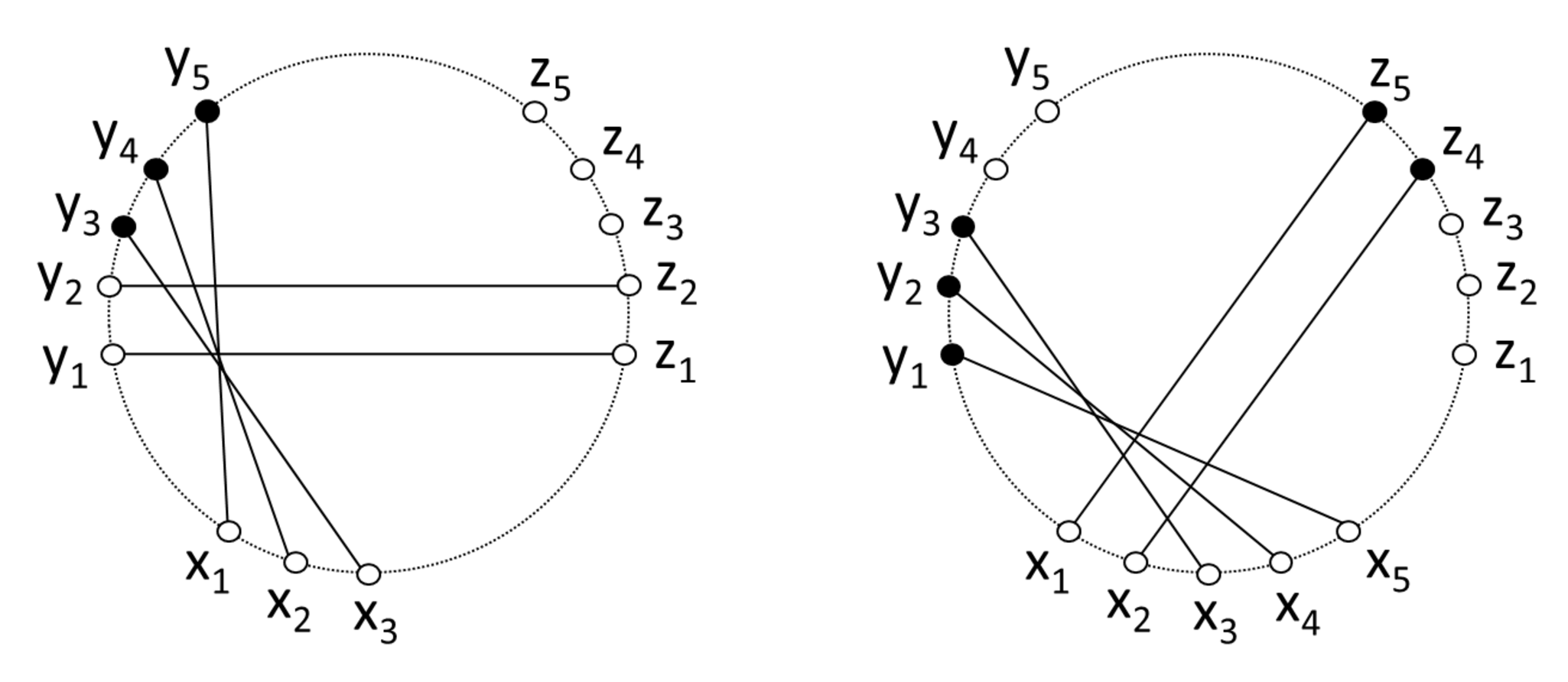}\\
\caption{$E((5,4,3), (0^3))$ (left) and $E((5,4,3,2,1), (1^2, 0^3))$ (right).}\label{fig2311_m_cross_n_parallel}
\end{center}
\end{figure}

\begin{figure}[H]
\begin{center}
\includegraphics[scale=0.3]{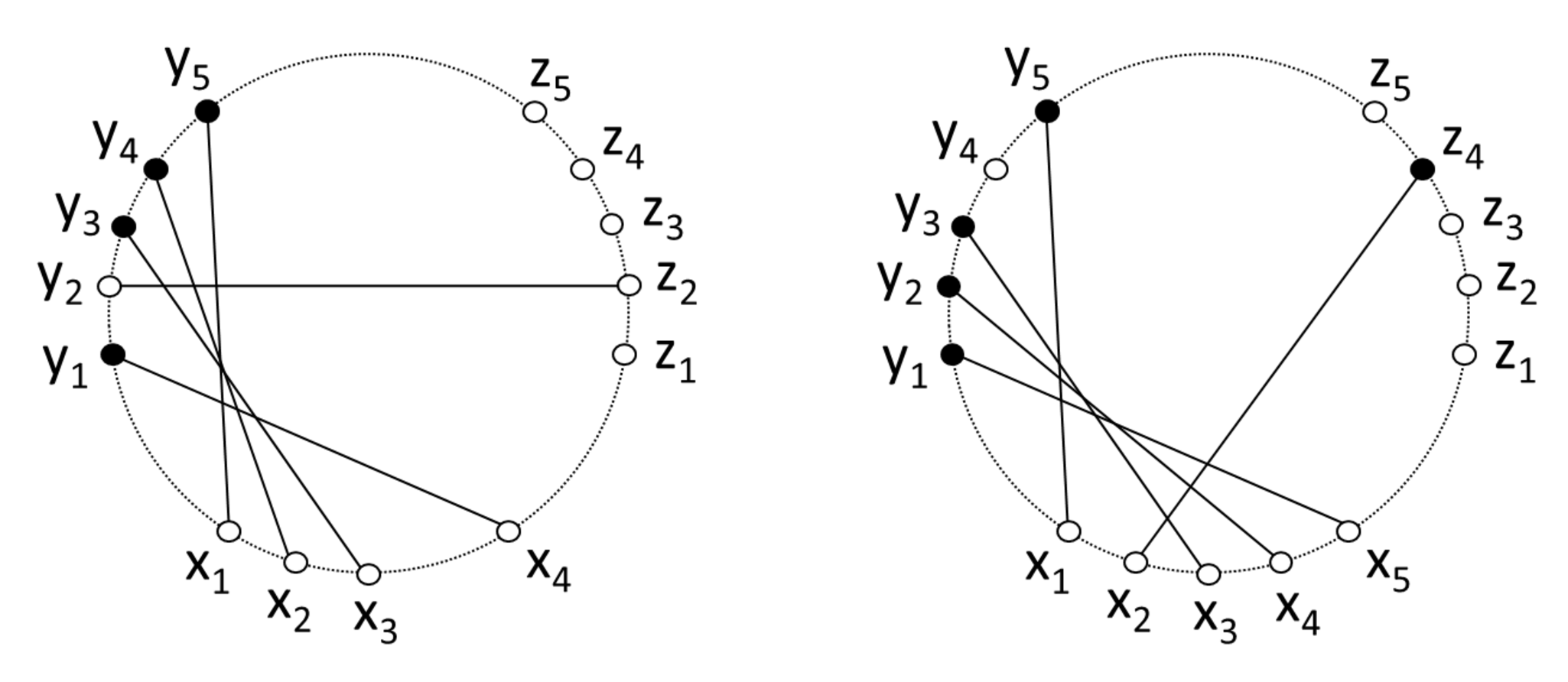}\\
\caption{$E((5,4,3,1), (0^4))$ (left) and $E((5,4,3,2,1), (0, 1, 0^3))$ (right).}\label{fig2301_E_n_k}
\end{center}
\end{figure}

Note that two different pairs $(\va_1, \vvarepsilon_1)$ and $(\va_2, \vvarepsilon_2)$ may give isomorphic chord diagrams.
For example, let $\va_1 = (n+m, n+m-1, \ldots, n+1)$, $\vvarepsilon_1 = (0^m)$ and
$\va_2 = (n+m, \ldots, 2, 1)$, $\vvarepsilon_2 = (1^n, 0^m)$.
Then, both ${\rm Int}(E(\va_1, \vvarepsilon_1))$ and ${\rm Int}(E(\va_2, \vvarepsilon_2))$ are isomorphic to $K_m + \overline{K_n}$, the join of $K_m$ and the compliment of $K_n$ (see Fig. \ref{fig2311_m_cross_n_parallel}).
Recall that $E_{n,k}$ is a chord diagram consisting of an $n$-crossing $C_n$ and an extra chord $e$, where $e$ crosses $k$ chords of $C_n$.
Let $\va_1 = (n+1, n, \ldots, n-k+2, n-k, n-k-1, \ldots, 1)$, $\vvarepsilon_1 = (0^n)$ and
$\va_2 = (n+1, n, \ldots, 1)$, $\vvarepsilon_2 = (0^{n-k}, 1, 0^k)$.
Then,  both $E(\va_1, \vvarepsilon_1)$ and $E(\va_2, \vvarepsilon_2)$ are also isomorphic to $E_{n,k}$ (see Fig. \ref{fig2301_E_n_k}).  

To exhibit the main results of the paper, we will introduce some more definitions. 
For a chord diagram $E$ and a vertex set $U$, $e_E(U)$ is defined as the number of chords of $E$ whose both endvertices are contained in $U$.
For $0 \le k \le n$, let $f(\va, \vvarepsilon, k)$ be defined as the number of nonintersecting chord diagrams $F \in {\cal NCD}(E(\va,\vvarepsilon))$ satisfying $k = e_F(Y) + e_F(X \cup Z)$, and
let $g(\va, \vvarepsilon, k)$ be defined as the number of $0$-$1$ Young diagrams $M$ of shape $\va$ satisfying the following conditions {\rm (1)}, {\rm (2)}, {\rm (3)}:\\
{\rm (1)~}Every column of $M$ contains at most one $1$.\\
{\rm (2)~}For $1 \le i \le m$, the $i$th row of $M$ contains $\varepsilon_i \pmod 2$ $1$'s.\\
{\rm (3)~}$M$ contains $k$ $1$'s.

Our main result is the following theorem.

\begin{thm}\label{thm_cdyd}
$f(\va, \vvarepsilon, k) = g(\va, \vvarepsilon, k)$.
\end{thm}

We have the following result from Theorem \ref{thm_cdyd} by taking the sum over all $k \ge 0$.

\begin{cor}\label{cor_cdyd}
$ex(E(\va,\vvarepsilon))$ equals the number of $0$-$1$ Young diagrams M of shape $\va$ satisfying the following conditions {\rm (1)}, {\rm (2)}:\\
{\rm (1)~}Every column of $M$ contains at most one $1$.\\
{\rm (2)~}For $1 \le i \le m$, the $i$th row of $M$ contains $\varepsilon_i \pmod 2$ $1$'s.
\end{cor}

By Corollary \ref{cor_cdyd}, we have the following result.

\begin{figure}[H]
\begin{center}
\includegraphics[scale=0.3]{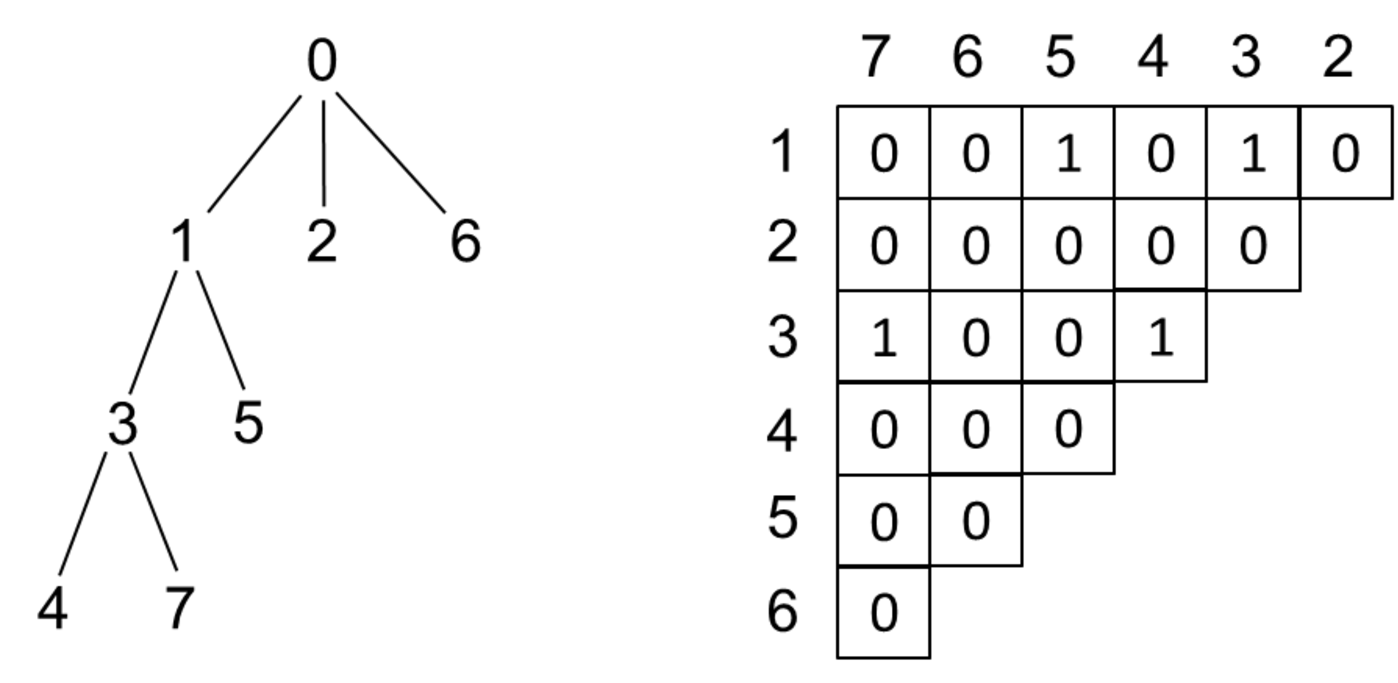}\\
\caption{An even increasing tree of order $7$(left) corresponds to a $0$-$1$ Young diagram of shape $(6, 5, 4, 3, 2, 1)$ satisfying the conditions of Corollary \ref{cor_n_crossing}.}\label{fig2301_even_increasing_tree}.
\end{center}
\end{figure}

\begin{cor}\label{cor_n_crossing}
$ex(C_n)$ equals the number of $0$-$1$ Young diagrams M of shape $(n, n-1, \ldots, 1 )$ satisfying the following conditions {\rm (1)}, {\rm (2)}:\\
{\rm (1)~}Every column of $M$ contains at most one $1$.\\
{\rm (2)~}For $1 \le i \le n$, the $i$th row of $M$ contains an even number of $1$'s.
\end{cor}
\noindent
{\bf Remark.} Corollary \ref{cor_n_crossing} could be derived differently by using even increasing trees, which were introduced by Kuznetsov, Pak, and Postnikov \cite{KPP1994}.
A rooted directed tree with the vertex set $\{ 0, 1, \ldots, n \}$ is called an {\it even increasing tree} of order $n+1$ if {\rm (1)} the root is $0$, {\rm (2)} for any given vertex $k \ne 0$, the labels of vertices are increasing along the path from $0$ to $k$, and {\rm (3)} for any given vertex $k \ne 0$, the outdegree of $k$ is even (see Fig. \ref{fig2301_even_increasing_tree}(left)).
 Kuznetsov, Pak, and Postnikov showed a bijection from the set of alternating permutations on $[n]$ to the family of even increasing trees of order $n+1$.
 On the other hand, there exists a bijection $\varphi$ from the family of even increasing trees of order $n+1$ to $0$-$1$ Young diagrams of shape $\va = (n, n-1, \ldots, 1)$ satisfying the conditions {\rm (1)} and {\rm (2)} in Corollary \ref{cor_n_crossing} as shown below.
 (The author owes the bijection to one of the anonymous reviewers.)
For a given even increasing tree $T$, if there is an edge $ij$ with $1 \le i < j \le n$, then we give $1$, otherwise, we give $0$, to the box of the $i$th row and the $(n+1-j)$th column of $\varphi(T)$ (see Fig. \ref{fig2301_even_increasing_tree}). 
Then $\varphi(T)$ is $0$-$1$ Young diagram of shape $\va$ satisfying the desired conditions.
It is not difficult to see that $\varphi^{-1}$ exists, hence $\varphi$ is a bijection.
\medskip

Again by Corollary \ref{cor_cdyd}, we have the following results.

\begin{cor}\label{cor_mn_crossing}
$ex(C_{m,n})$ equals the number of $0$-$1$ Young diagrams M of shape $ (n+m, n+m-1, \ldots, n+1)$ satisfying the following conditions {\rm (1)}, {\rm (2)}:\\
{\rm (1)~}Every column of $M$ contains at most one $1$.\\
{\rm (2)~}For $1 \le i \le m$, the $i$th row of $M$ contains an odd number of $1$'s.
\end{cor}

\begin{cor}\label{cor_E_n_k}
For $1 \le k \le n$, $ex(E_{n,k})$ equals the number of $0$-$1$ Young diagrams M of shape $(n+1, n, \ldots, n-k+2, n-k, n-k-1, \ldots, 1)$ satisfying the following conditions {\rm (1)}, {\rm (2)}:\\
{\rm (1)~}Every column of $M$ contains at most one $1$.\\
{\rm (2)~}For $1 \le i \le n$, the $i$th row of $M$ contains an even number of $1$'s.
\end{cor}
\noindent
{\bf Remark.} It is shown that some types of chord diagrams $E(\va, \vvarepsilon)$ are corresponding to some families of $0$-$1$ Young diagrams through chord expansions.
It may be interesting to find a more general correspondence between an arbitrary chord diagram and some more general objects containing or similar to $0$-$1$ Young diagrams.
\medskip

We remark an identity concerning the chord expansion number as an application of Corollary \ref{cor_cdyd}.

\begin{cor}\label{sum_identity}
Let $\va = (a_1, a_2, \ldots, a_m)$ with $a_1 > a_2 > \cdots > a_m \ge 1$.
Let $a_{m+1} = 0$.
Then $\sum_{\vvarepsilon \in \{ 0,1 \}^m} ex(E(\va,\vvarepsilon)) = \prod_{i=1}^{m} (i+1)^{a_i - a_{i+1}} $.
\end{cor}
{\bf Proof. }
By Corollary \ref{cor_cdyd}, taking the sum over all $\vvarepsilon \in \{ 0,1 \}^m$, the left-hand side of the equation of the claim is the number of all $0$-$1$ Young diagrams $M$ of shape $\va$ in which each column has at most one $1$.
This is counted by the product of $c_j + 1$ for $1 \le j \le n$, where $c_j$ is the length of the $j$th column of $M$, which is $\prod_{j=1}^{n} (c_j + 1)$.
Since the number of indices $j$ satisfying $c_j = i$ is $a_i - a_{i+1}$ for $0 \le i \le m$, we have $\prod_{j=1}^{n} (c_j + 1) = \prod_{i=1}^{m} (i+1)^{a_i - a_{i+1}}$.  
\owari
\medskip

\section{Proof of Theorem \ref{thm_cdyd}}
Let us define a total order $\prec$ for the set of monotonically decreasing finite sequences of positive integers.
For $\va = (a_1, a_2, \ldots, a_s)$ and $\vb = (b_1, b_2, \ldots, b_t)$, $\va \prec \vb$ if and only if {\rm (1)} $s < t$ or {\rm (2)} $s = t$, $a_i = b_i$ for all $1 \le i < r$ and $a_r < b_r$ for some $r$.
We proceed by induction on the order $\prec$.
Firstly, assume that $\va = (1)$, the minimum element of the order $\prec$, and $\vvarepsilon = (\varepsilon_1)$.
In this case, $F = E(\va, \vvarepsilon)$ consists of one chord and $e_F(Y) + e_F(X \cup Z) = \varepsilon_1$.
Hence, we have $f(\va, \vvarepsilon, k) = \delta_k^{\varepsilon_1}$ for $k=0, 1$, where $\delta_i^j$ is Kronecker delta.
On the other hand, the only $0$-$1$ Young diagram of shape $\va$ satisfying the conditions {\rm (1)}, {\rm (2)}, {\rm (3)} has one box which contains $\varepsilon_1$ $1$'s.
Hence, we have $g(\va, \vvarepsilon, k) = \delta_k^{\varepsilon_1}$ for $k=0, 1$, as claimed.

\begin{figure}[H]
\begin{center}
\includegraphics[scale=0.3]{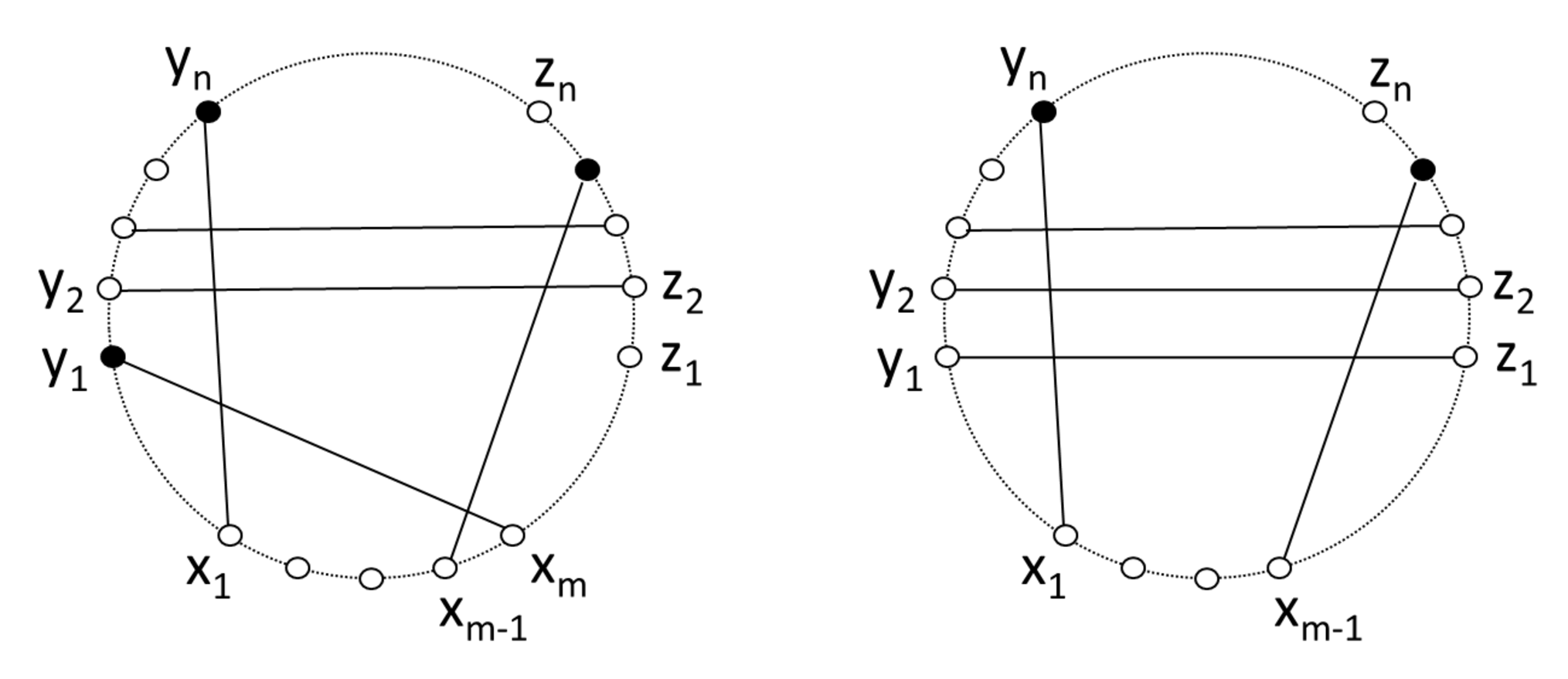}\\
\caption{$E(\va, \vvarepsilon)$ (left) and $E(\va', \vvarepsilon')$ (right) in Case $1$. $1$.}\label{fig2311_thm1_case11}
\end{center}
\end{figure}

Let $\va = (a_1, a_2, \ldots, a_m)$ with $(1) \prec \va$.
\medskip\\
\underline{Case $1$.} $a_m = 1$.
\medskip\\
\underline{Case $1$. $1$.} $\varepsilon_m = 0$.
\medskip\\
In this case, we have $x_m y_1 \in E(\va, \vvarepsilon)$.
By replacing $x_m y_1$ with $y_1 z_1$, we have $E(\va', \vvarepsilon')$, where $\va' = (a_1, \ldots, a_{m-1})$ and $\vvarepsilon' = (\varepsilon_1, \ldots, \varepsilon_{m-1})$ (see Fig. \ref{fig2311_thm1_case11}).

\begin{figure}[H]
\begin{center}
\includegraphics[scale=0.3]{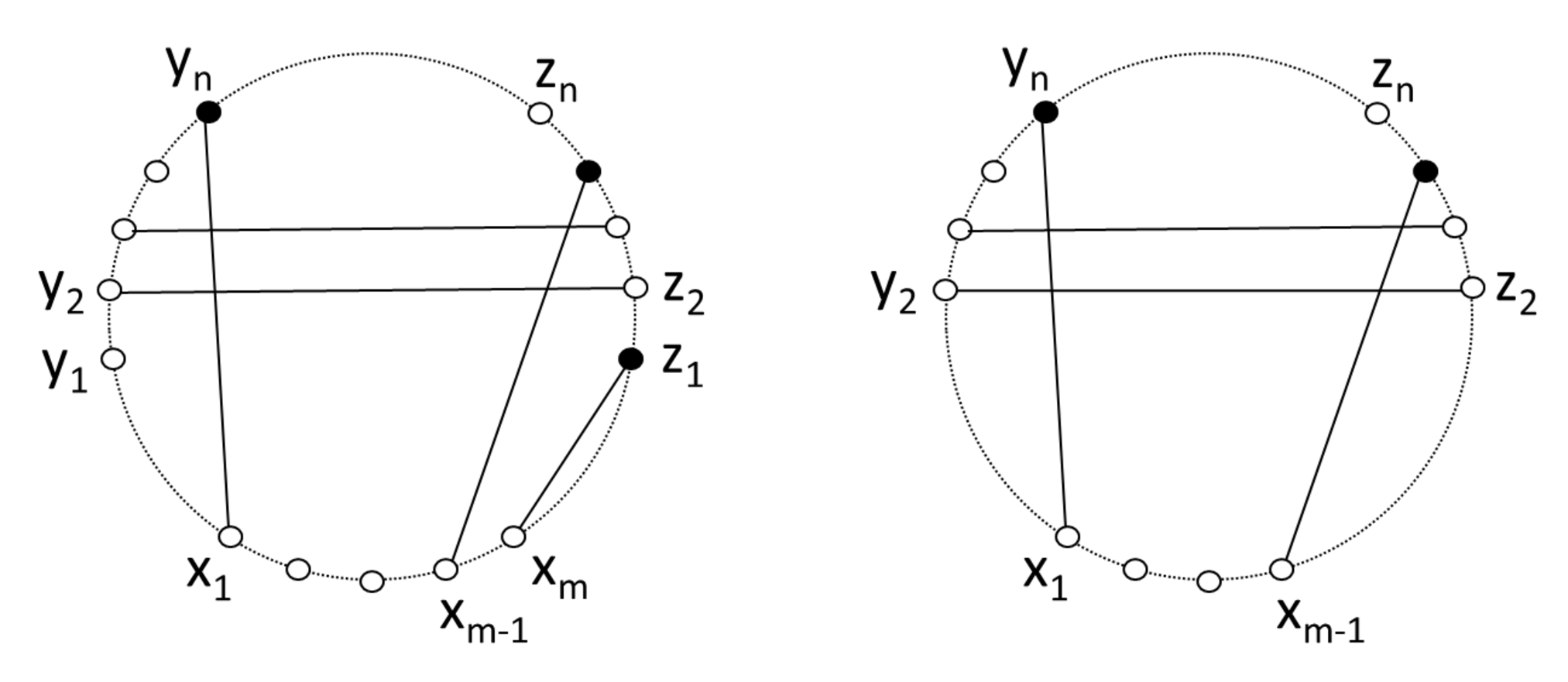}\\
\caption{$E(\va, \vvarepsilon)$ (left) and $E(\va', \vvarepsilon')$ (right) in Case $1$. $2$.}\label{fig2311_thm1_case12}
\end{center}
\end{figure}

Since $E(\va,\vvarepsilon)$ and $E(\va', \vvarepsilon')$ are isomorphic and both $x_m$ and $z_1$ are contained in $X \cup Z$, 
we have $f(\va, \vvarepsilon, k) = f(\va', \vvarepsilon', k)$.
On the other hand, for a $0$-$1$ Young diagram of shape $\va$, its $m$th row has one box which contains $0$.
Hence, we have $g(\va, \vvarepsilon, k) = g(\va', \vvarepsilon', k)$. 
By inductive hypothesis, we have $f(\va, \vvarepsilon, k) = f(\va', \vvarepsilon', k) = g(\va', \vvarepsilon', k) = g(\va, \vvarepsilon, k)$. 
\medskip\\
\underline{Case $1$. $2$.} $\varepsilon_m = 1$.
\medskip\\
In this case, we have $x_m z_1 \in E(\va,\vvarepsilon)$ and $x_m z_1$ is an isolated chord.
By removing $x_m z_1$ from $E(\va,\vvarepsilon)$, we have $E(\va',\vvarepsilon')$, where $\va' = (a_1-1, a_2-1, \ldots, a_{m-1}-1)$ and $\vvarepsilon' = (\varepsilon_1, \ldots, \varepsilon_{m-1})$ (see Fig. \ref{fig2311_thm1_case12}).

\begin{figure}[H]
\begin{center}
\includegraphics[scale=0.3]{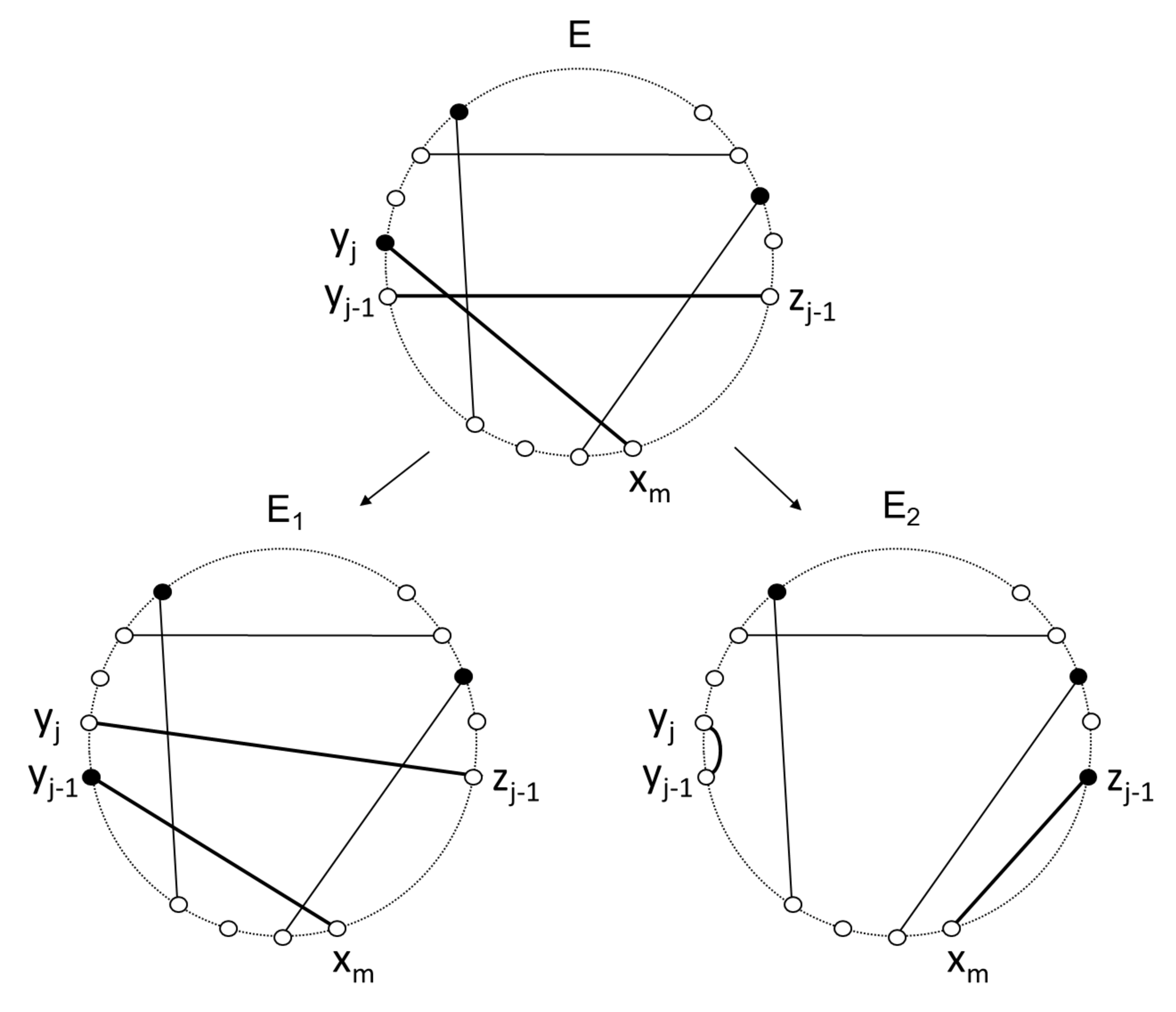}\\
\caption{$E = E(\va, \vvarepsilon)$ (top) , $E_1$ (left) and $E_2$ (right) in Case $2$. $1$.}\label{fig2311_thm1_case211}
\end{center}
\end{figure}

Since $x_m z_1$ is an isolated chord and both endvertices of $x_m z_1$ are contained in $X \cup Z$, we have $f(\va, \vvarepsilon, k) = f(\va', \vvarepsilon', k-1)$.
On the other hand, for a $0$-$1$ Young diagram of shape $\va$ with $\varepsilon_m = 1$, its $m$th row has one box which contains $1$ and the $(i, 1)$ elements are $0$ for $1 \le i \le m-1$.
Hence, we have $g(\va, \vvarepsilon, k) = g(\va', \vvarepsilon', k-1)$.
By inductive hypothesis, we have $f(\va, \vvarepsilon, k) = f(\va', \vvarepsilon', k-1) = g(\va', \vvarepsilon', k-1) = g(\va, \vvarepsilon, k)$.
\medskip\\
\underline{Case $2$.} $a_m \ge 2$.
\medskip\\
Put $j = a_m$.
\medskip\\
\underline{Case $2$. $1$.} $\varepsilon_m = 0$.
\medskip\\
In this case, $x_m y_j$ and $y_{j-1} z_{j-1}$ are contained in $E(\va, \vvarepsilon)$.
Let us consider a chord expansion of $E(\va, \vvarepsilon)$ with respect to $S = \{ x_m y_j, y_{j-1} z_{j-1} \}$.
The resultant chord diagrams are $E_1 = (E(\va, \vvarepsilon) \setminus S) \cup \{ y_j z_{j-1}, x_m y_{j-1} \}$ and $E_2 = (E(\va, \vvarepsilon) \setminus S) \cup \{ y_j y_{j-1}, x_m z_{j-1} \}$  (see Fig. \ref{fig2311_thm1_case211}).

\begin{figure}[H]
\begin{center}
\includegraphics[scale=0.3]{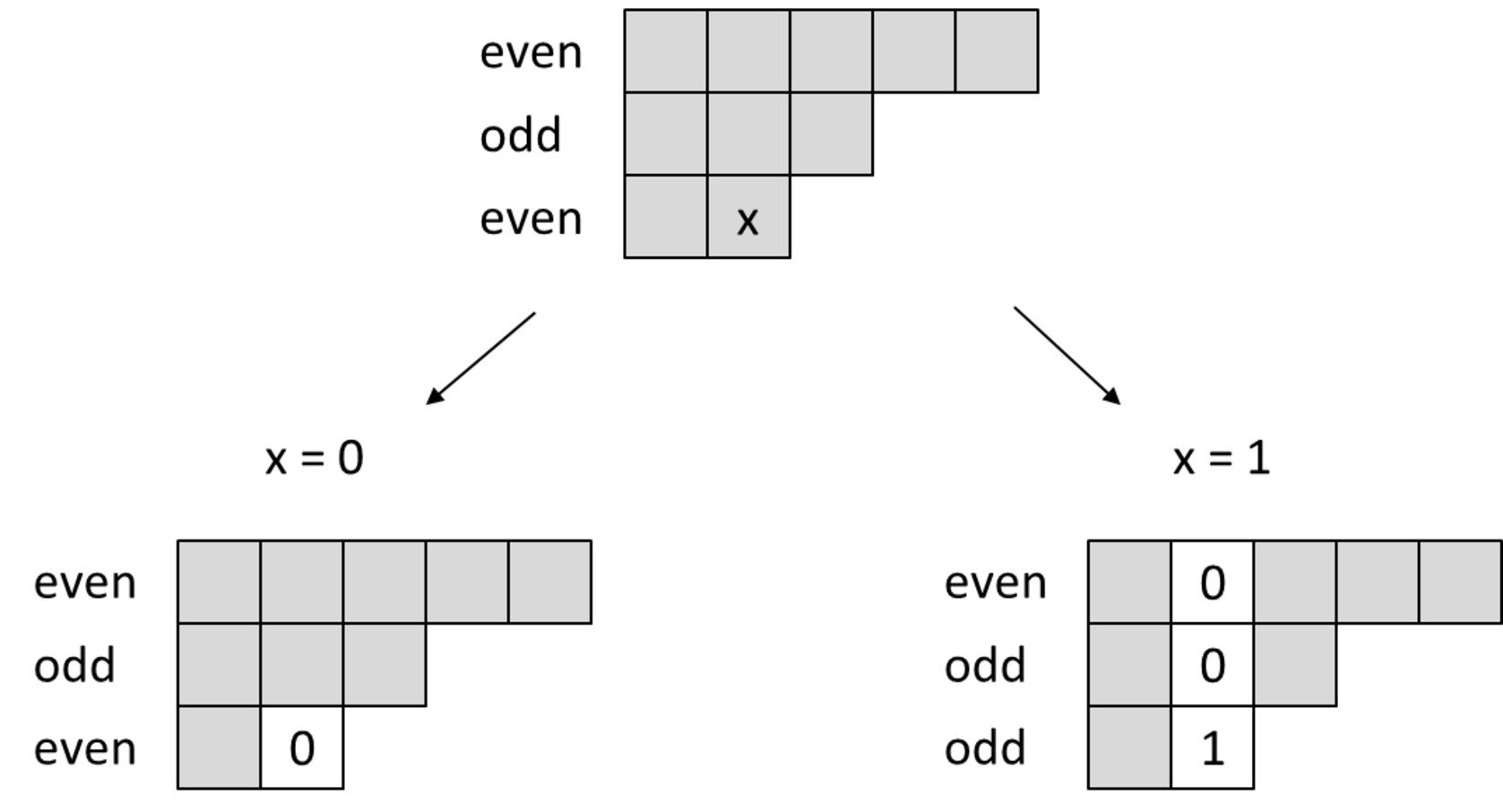}\\
\caption{Sets of $0$-$1$ Young diagrams (illustrated by shaded boxes) corresponding to chord diagrams of Fig. \ref{fig2311_thm1_case211}.}\label{fig2301_thm1_case212}
\end{center}
\end{figure}

Then we have $E_1$ is isomorphic to $E(\va_1, \vvarepsilon_1) $, where $\va_1 = (a_1, \ldots, a_{m-1}, a_m-1)$ and $\vvarepsilon_1  = \vvarepsilon$.
Furthermore, $E_2$ is isomorphic to $E(\va_2, \vvarepsilon_2) \cup \{ y_j y_{j-1} \}$, where $\va_2 = (a_1-1, a_2-1, \ldots, a_m-1)$ and $\vvarepsilon_2  =  (\varepsilon_1, \ldots, \varepsilon_{m-1}, 1-\varepsilon_m)$.
Note that $y_j y_{j-1}$ is an isolated chord, whose both endvertices are contained in $Y$.
Hence, we have $f(\va, \vvarepsilon, k) =  f(\va_1, \vvarepsilon_1, k) + f(\va_2, \vvarepsilon_2, k-1)$.
On the other hand, for a $0$-$1$ Young diagram of shape $\va$, let us consider the value $x$ of the $(m, j)$ element (see Fig. \ref{fig2301_thm1_case212}).
If $x=0$, then the boxes in the $j$-column with entries from $1$ to $m-1$ have at most one $1$; also the boxes in the $m$-row with entries from $1$ to $j-1$ have an even number of $1$'s.
If $x=1$, then the boxes in the $j$-column with entries from $1$ to $m-1$ have $0$; also the boxes in the $m$-row with entries from $1$ to $j-1$ have an odd number of $1$'s.
Hence, we have $g(\va, \vvarepsilon, k) = g(\va_1, \vvarepsilon_1, k) + g(\va_2, \vvarepsilon_2, k-1)$.   
Therefore, by the inductive hypothesis, we have 
\begin{eqnarray*}
f(\va, \vvarepsilon, k) & = & f(\va_1, \vvarepsilon_1, k) + f(\va_2, \vvarepsilon_2, k-1)\\
 & = & g(\va_1, \vvarepsilon_1, k) + g(\va_2, \vvarepsilon_2, k-1)\\
 & = & g(\va, \vvarepsilon, k).
\end{eqnarray*}
\medskip\\
\underline{Case $2$. $2$.} $\varepsilon_m = 1$.
\medskip\\
In this case, $x_m z_j$ and $y_{j-1} z_{j-1}$ are contained in $E(\va, \vvarepsilon)$.
In the same way as in Case 2. 1, by exchanging the roles of $y_{j-1}, y_j$ and $z_{j-1}, z_j$, respectively, 
 it is proved that $f(\va, \vvarepsilon, k) = g(\va, \vvarepsilon, k)$.  
\owari
\medskip

{\bf Acknowledgment.} The author would like to thank anonymous reviewers for their valuable comments, which greatly improved the earlier version of this paper.
     

\end{document}